\documentclass[12pt]{article}
\usepackage{amsmath,amsfonts,amssymb,amsthm}
\usepackage{ae,aecompl}

\begin{document}
\baselineskip18pt

\newcommand{\A}{\mathcal{A}}        
\renewcommand{\H}{\mathcal{H}}      
\renewcommand{\t}{\tau}      
\newcommand{\Tb}{\mathbb{T}}         
\newcommand{\bZ}{\mathbb{Z}}         
\newcommand{\bN}{\mathbb{N}}         
\newcommand{\bR}{\mathbb{R}}         
\newcommand{\bC}{\mathbb{C}}         
\newcommand{\te}{\theta}            
\newcommand{\cO}{\mathcal{O}}        
\newcommand{\Dslash}{{D\mkern-11.5mu/\,}} 
\def\vol{{\mathrm{vol}}}


\begin{center}

{\Large
{\bf
SPINORS AND THETA DEFORMATIONS
}
}
~\\
{\Large
~\\
{\bf Ludwik D\k{a}browski}\\
~\\
{\it
SISSA, Via Beirut 2-4, I-34014, Trieste, Italy
}
}
\end{center}
~\\
\begin{abstract}
\baselineskip17pt
\textwidth 8cm
\noindent
The construction \cite{cl01} of Dirac operators on $\theta$-deformed manifolds 
is recalled, stressing the aspect of spin structure.
The description of \cite{CD-V} is extended to arbitrary spin structure.\\
\end{abstract}

\vspace{1pc}

\newpage

%

In noncommutative geometry \`a la Connes,
{\em (compact) Riemannian, spin geometry} is described using a
{\em (unital) spectral triple} $(\A,\H,D)$, where $\A$ is a $*$-algebra, 
$\H$ is a Hilbert space  
and $D$ is 
an operator on $\H$  such that 
$D\!  =\!  D^\dagger$, ~$(D\! -\! z)^{-1}$ is compact and $[D, a]$ is bounded 
$\forall a\in A$.
Classically, for a manifold $M$ with Riemannian 
metric $g$ and spin structure $\sigma$,
the canonical spectral triple is given by
$\A=C^\infty (M, \mathbb C)$,  
$\H=L^2 (\Sigma_\sigma, vol_g)$ and $D$ being the Dirac operator $\Dslash$ 
constructed from the Levi-Civita connection (metric-preserving and torsion-free). 

Using further seven conditions, Connes stated his reconstruction conjecture
(theorem if $\A= C^\infty (M)$, see also [RV]).
These conditions 
({\it regularity} (or {\it smoothness}), {\it dimension} (or {\it summability}),
{\it finiteness and absolute continuity}, {\it reality}, {\it first order},
{\it orientability}, {\it Poincar\'e duality}) have 
been formulated [C] also for noncommutative manifolds.
They are satisfied e.g. by finite geometries, the noncommutative torus and more 
generally by the isospectral $\theta$-deformations of the canonical spectral triples
on a Riemannian spin manifold $M$.

For the latter example one assumes that the metric $g$ on $M$ is such that its 
isometry group contains the two-dimensional torus, ${\rm Isom}(M,g)\supset \Tb^2 $
(so one considers $\Tb^2$-invariant metrics only).
Any $\t = (\t_1 ,\t_2)\in \Tb^2$, $\t_1=e^{is_1}, \t_2 = e^{is_2}$, 
lifts to a pair of automorphisms 
on $\H = L^2(\Sigma )$,
which altogether form a double covering 
$$\widetilde\Tb^2\ni {\tilde \t} \mapsto \t \in \Tb^2 ,$$
the kernel of which we denote by ${\bZ'_2}$.
This covering must be one of the exactly four (inequivalent) double coverings 
$\widetilde\Tb^2_{j_1,j_2}\to\Tb^2 $, numbered by $j_1,j_2=0,1$:\\
\begin{equation}\label{cover}
\begin{array}{c}
\Tb^2 \times \bZ_2 ~\ni
~ (\t_1,\t_2,\pm 1)\mapsto (\t_1,\t_2) ~~{\rm if}~ j_1,j_2=0,0 \ ,  \\
\\
\Tb^2~ \ni
~ (\t_1,\t_2)\mapsto ((\t_1)^2,\t_2) ~~{\rm if}~ j_1,j_2=1,0 \ , \\
\\
\Tb^2 ~\ni
~ (\t_1,\t_2)\mapsto (\t_1,(\t_2)^2) ~~{\rm if}~ j_1,j_2=0,1\ ,  \\
\\
\Tb^2 /\bZ_2^{diag} ~\ni
~ [\t_1,\t_2] \mapsto ((\t_1)^2,(\t_2)^2) ~~{\rm if}~ j_1,j_2=1,1\ ,
\end{array}\end{equation}
where $[\t_1,\t_2]=\{(\t_1,\t_2), (-\t_1,-\t_2)\}$.
Note that only three different covering manifolds occur here (since $\widetilde\Tb^2_{0,1}$ and $\widetilde\Tb^2_{1,0}$ are identical) -
we shall often call them untwisted, twisted along one loop, or twisted  
along two loops, respectively.
However the four covering maps are all different and their winding numbers
along the two factors $S^1$ of $\Tb^2$ are $j_1+1$ and $j_2+1$, 
respectively. 

In what follows we shall need the action (of the generator) of the kernel ${\bZ'_2}$
on $\widetilde\Tb^2_{j_1,j_2}$
\begin{equation}\label{Z2prime}
\begin{array}{c}
 (\t_1,\t_2,\pm 1)\mapsto (\t_1,\t_2,\mp 1) ~~{\rm if}~ j_1,j_2=0,0 \ ,  \\
\\
 (\t_1,\t_2)\mapsto (-\t_1,\t_2) ~~{\rm if}~ j_1,j_2=1,0 \ , \\
\\
 (\t_1,\t_2)\mapsto (\t_1,-\t_2) ~~{\rm if}~ j_1,j_2=0,1\ ,  \\
\\
~ [\t_1, \t_2 ] \mapsto [-\t_1,\t_2] ~~{\rm if}~ j_1,j_2=1,1\ .
\end{array}\end{equation}
Which one of these coverings occurs depends on $M$ and the spin structure $\sigma$. 
Each of them is realized for some $M$ and $\sigma$.\\
Example: On $M= S^N$, $N\geq 3$, there is a unique  $\sigma$, 
for which $\widetilde\Tb^2 = \widetilde\Tb^2_{1,1}$.\\
Example: On $M\! =\! U(2) $, there are two such $\sigma$, 
and $\widetilde\Tb^2$ is $\widetilde\Tb^2_{01}$ or $\widetilde\Tb^2_{10}$.\\
Example: On $M= \Tb^N $, $N\geq 2$, there are $2^N$ such ~$\sigma$, 
but the case of $\Tb^2$-action explores only $2^2 $ of them
(the action of $\Tb^N $ can be used though, to be discussed later on)
so the situation can be understood already for $N=2$.\\
\underline{Example}: On $M= \Tb^2$, there are $2^2=4$ spin structures 
$\sigma_{j_1,j_2}$, $j_1,j_2=0,1$. 
They are given explicitly by (trivial) principal bundles 
together with four morphisms as follows
$$
\def\normalbaselines{\baselineskip30pt
\lineskip3pt \lineskiplimit3pt }
\def\mapright#1{\smash{
\mathop{\!\!\!-\!\!\!\longrightarrow\!\!\!}
\limits^{#1}}}
\def\mapup#1{\Big\uparrow\rlap{$\vcenter{\hbox{$\scriptstyle#1$}}$}}
\def\mapdn#1{\Big\downarrow\rlap{$\vcenter{\hbox{$\scriptstyle#1$}}$}}
\def\mapsw#1{\smash{\mathop{\vector(-1,-1){20}}\limits^{#1}}}
\def\mapse#1{\smash{\mathop{\vector(1,-1){20}}\limits^{#1}}}
  \begin{array}{ccc}
 \Tb^2 \times Spin(2) & \ni &\!\!\!\!\!\!\! (e^{ix_1},e^{ix_2},e^{i\phi}) \cr
& & \cr
 \setlength{\unitlength}{1mm}
\begin{picture}(12,0)
\put(2,5){\vector(0,-1){8}}
\put(4,0){\mbox{\scriptsize$ \eta_{j_1,j_2}$}}
\end{picture}
&  &
 \setlength{\unitlength}{1mm}
\begin{picture}(12,0)
\put(2,5){\vector(0,-1){8}}
\put(4,0){\mbox{\scriptsize$ \eta_{j_1,j_2}$}}
\end{picture}
 \cr
& & \cr
 \Tb^2 \times SO(2)&\ni  &  (e^{ix_1},e^{ix_2},R(j_1x_1+j_2x_2+2\phi)). \cr
\end{array}
%
$$
Here, $Spin(2) = U(1)$ and $R(r)\in SO(2)$ is the rotation by the angle $r$.
It is easy to see that the tautological (isometric) action of 
$\t = (e^{is_1}, e^{is_2}) \in \Tb^2$,
$$\alpha_{\t}: (e^{ix_1},e^{ix_2})\mapsto (e^{i(s_1+x_1)},e^{i(s_2+x_2)})$$
on $M=\Tb^2$ (equipped with the usual flat metric)
has exactly two lifts
\begin{equation}\label{lift}
(e^{ix_1},e^{ix_2},e^{i\phi})\mapsto (e^{i(s_1+x_1)},e^{i(s_2+x_2)},
\pm e^{-i(j_1s_1+j_2s_2)/2}e^{i\phi}),
\end{equation}
which form a subgroup of automorphisms of the (trivial) principal bundle 
$ \Tb^2 \times Spin(2)$. 
The automorphisms of $\H$ associated with them
(via pullback if $\psi\in\H$ is regarded as an equivariant function on 
$ \Tb^2 \times Spin(2)$,
but we can regard $\psi$ just as functions on $ \Tb^2$)
form $\widetilde\Tb^2_{j_1,j_2}$.
As far as the dependence on  $(e^{is_1}, e^{is_2})$ is regarded this
determines the two (local) sections of the coverings \eqref{cover}.
Varying $j_1,j_2=0,1$ yields exactly all the four coverings  \eqref{cover}.\\

Now we summarize the theta  (isospectral) deformation  
\cite[Section~V]{cl01} of the canonical spectral triple. 
Denote by $ \widehat\alpha_{\tilde \t} $ the action of the (unitary) operator on $\H$ 
corresponding to $\tilde\t$.
It satisfies
$$   
\widehat\alpha_{\tilde \t}(a\psi)=\alpha^*_{\t}(a) \widehat\alpha_{\tilde \t}(\psi), \forall a\in \A,  \psi\in\H
$$
and
$$
[\widehat\alpha_{\tilde \t}, D] = 0, \, \, [\widehat\alpha_{\tilde \t}, J] = 0.
$$
This induces a bigrading of $B(\H)$, according to which $K\in B(\H)$ has bidegree $(n_1,n_2)$ if
\begin{equation}
\widehat\alpha_{\tilde s}\, K \, \widehat\alpha_{\tilde s}^{-1} = 
\exp(is_1 n_1 + is_2 n_2) \, K \, , \, \, \, \forall \, {\tilde s} 
\in \widetilde\Tb^2 \, .
\end{equation}
Any smooth $K$, i.e. \!\!\! such that 
$s \rightarrow \widehat\alpha_{\tilde s} \, K \, \widehat\alpha_{\tilde s}^{-1}$ is $C^\infty$ in $\| \,\|$,
has a unique norm-convergent decomposition 
\begin{equation}
K = \sum_{n_1,n_2} \, 
{K}_{n_1,n_2} \, ,
\end{equation}
with ${K}_{n_1,n_2}$ of bidegree $(n_1,n_2)$ and 
$||{K}_{n_1,n_2} ||$ sequence of rapid decay in $(n_1,n_2)$,  and conversely.
Write (locally)
%
\begin{equation}
\widehat\alpha_{\tilde s} = \exp(is_1 p_1 + is_2 p_2) \, 
\end{equation}
where the generators $p_\ell$ satisfy\\
\centerline{$D p_\ell =p_\ell D $, $J p_\ell = -p_\ell J$,}
\centerline{${\rm Spec}(p_1 )\in \bZ + j_1/2 $, 
${\rm Spec}(p_2 )\in \bZ + j_2/2 $ if $\widetilde\Tb^2 = \widetilde\Tb^2_{j_1,j_2}$.}\\
~\\
Let $\lambda = \exp(2 \pi i \theta)$. For smooth $K$, define   $\widehat K $ by 
\begin{equation}
\widehat 
K = \sum_{n_1,n_2} \, {K}_{n_1,n_2} \, \lambda^{n_2 p_1} 
\end{equation}
(which is still convergent). Then
$$
\widehat 
K_{n_1,n_2} \widehat 
K'_{n_1',n_2'} 
=
{\widehat{(K_{n_1,n_2}*K'_{n_1',n_2'})}}\ ,
$$
where
$$
{K_{n_1,n_2}}*{K'_{n_1',n_2'}} 
=
\lambda^{n_1' n_2} {K_{n_1,n_2} \, K'_{n_1',n_2'}}. 
$$
The algebra $\widehat \A = : C^\infty(M_\te) \, (=\A_\te) $ 
is called $\theta$-deformation of $\A$.
Its continuous version $C(M_\te)$
is nothing but the Rieffel 
deformation quantization \cite{R} of $C(M)$.
Clearly this gives  $C^\infty(M_0) = C^\infty(M)$ when $\theta = 0$.
Moreover $C^\infty((\Tb^2)_\te) = C^\infty(\Tb^2_\te)$ 
(smooth algebra of the noncommutative torus),
as it should.\\
\\
\noindent
Next, if $J$ is the canonical real structure (charge conjugation) for $\A$, then
\begin{equation}
\widehat{J} := J \, \lambda^{- p_1 p_2} = \lambda^{p_1 p_2} \, J 
\,  \label{realtwist}
\end{equation}
is a real structure for $\widehat \A$.
Moreover if dim$M$ is even and $M$ is oriented, there is a grading 
$$\widehat\gamma=\gamma .$$
  
The following describes the theta deformation 
of the canonical spectral triple.\\
~\\
\noindent
{\bf Proposition}~[CL]. The datum $(C^\infty(M_\te), \H , D)$ 
with the real structure $\widehat{J}$ (and $\gamma$, if any)
is a spectral triple which satisfies the seven conditions 
required of a noncommutative manifold.\\

In another interesting paper \cite{CD-V} a multiparameter deformation $S^{N-1}_\theta$ of $S^{N-1}$ and $\bR^N_\theta$ of $\bR^N$ were introduced,
where $\theta$ is a $N\times N$ real skew-symmetric matrix of parameters.
Moreover they were shown to 
have a more `functorial' realization. 
Namely, when $M=S^{N-1}$ or  $\bR^N$, $C^\infty(M_{\theta})$ admits an action $\alpha^*$ of $\Tb^N$ and there is a `splitting'  
isomorphism 
\begin{equation}\label{splitiso}
\kappa:C^\infty(M_{\theta}) \approx 
\left(
C^\infty(M)\widehat\otimes C^\infty(\Tb^N_{\theta})
\right)^{ \alpha^*\otimes \beta^{-1}},
\end{equation}
with the fixed point subalgebra of the 
action $\alpha^*\otimes\beta^{-1}$ of $\Tb^N$,
given on the generators by 
\begin{equation}\label{split}
\kappa (z_{m,\te})={z_m}\otimes	u_{m,\te}, \,\forall m ,
\end{equation}
where $u_{m,\te}$ is the $m$-th generator of $\Tb^N_{\theta}$.
%
Here $\widehat\otimes$ denotes the (unique) completion of $\otimes$
and $\beta $  denotes the usual action of $\Tb^N$ on $C^\infty(\Tb^N_{\theta})$.\\
This holds not just for $S^N$ and $\bR^N$ but in fact
for any $M$ such that ${\rm Isom}(M,g)\supset \Tb^N$ [CD-V, sect. 12].
Namely, there still exists a splitting isomorphism $\kappa$ of \eqref{splitiso},
which using the fact that $C^\infty(M_{\theta})$ is $\bZ^N$-graded,
we write explicitly as 
\begin{equation}\label{splitt}
\kappa (a_{\delta_{m,\te}})=a_{\delta_{m}} \otimes	u_{m,\te}
\end{equation}
on the elements with the $N$-grade 
$\delta_{m}$, where $(\delta_{m})_k= 1$
if $m=k$ and $0$ if $m\neq k$
(and extend it to $C^\infty(M_{\theta})$).
The formula \eqref{split} is recovered by noting that $z_{m,\te}$ has grade $\delta_{m}$.
Thus symbolically, as `virtual spaces', the $\theta$-deformation corresponds to
gluing the noncommutative torus along the action of the torus, 
$$M_{\theta} = (M\times \Tb^N_{\theta})/\Tb^N.$$
%

Next, [CD-V] generalize 
\eqref{splitiso} to (smooth sections of) $\Tb^N$-equivariant vector bundles.
Let $S$ be a  smooth vector bundle over $M$,
which admits a lift $\widehat\alpha$ of isometries $\Tb^N$ of $M$, so
for any $\t\in\Tb^N$ there is 
$\widehat\alpha_\t$ such that
     $$\widehat\alpha_{\t}(a\psi)=\alpha_{\t}(a) \widehat\alpha_{\t}(\psi),
     \quad\forall a\in C^\infty(M), \psi\in C^\infty(M,S)$$
(this happens for `tensorial' bundles). 
Then
$$
C^\infty(M_{\theta},S) := 
\left(C^\infty(M,S)\widehat\otimes C^\infty(\Tb^N_{\theta})\right)^{\widehat\alpha\otimes \beta^{-1}}
$$
is a (finite projective) diagonal bimodule over $C^\infty(M_{\theta})$.
\noindent 
Moreover any operator $K$ on $C^\infty(M,S)$,
such that  $[K, \widehat\alpha_{s}]=0$, defines an operator $K_{\theta}$ on  
$C^\infty(M_{\theta},S)$ by
$$K_{\theta}:=K\otimes I\downharpoonright C^\infty(M_{\theta},S).$$
If $K$ is an order $n$ differential operator, $K_{\theta}$ is an order $n$ 
differential operator on the bimodule $C^\infty(M_{\theta},S)$ over 
$C^\infty(M_{\theta})$.
(This was employed for $S=\wedge \Tb^\ast M$ and $K=d$ and for the Hodge star 
to show that the Hochschild dimension and Poincar\`e duality are maintained under 
$\te$-deformations).

Now, as said before, concerning `spinorial' vector bundles such as $\Sigma$, 
there is a subtlety related to the fact that they carry the action 
$ \widehat\alpha_{\tilde \t}$
of a twofold covering $\tilde \Tb^N$ of $\Tb^N$, ${\tilde \t}\mapsto \t$, 
such that ~$ \widehat\alpha_{\tilde\t}(a\psi)=\alpha^*_{\t}(a) \widehat\alpha_{\tilde\t}(\psi) $,
rather than of $\Tb^N$.
To proceed further one thus needs to know 
the deformation $(\tilde \Tb^N)_\te$.
The answer provided in \cite{CD-V} is
\begin{equation}
\label{spadesuit}
\tilde \Tb^N_\te  \equiv  (\tilde \Tb^N)_\te := \Tb^N_{\frac{1}{2}\theta}.
\end{equation}
(Its status will be explained later on).
According to \cite{CD-V}, $C^\infty(\Tb^N_{\frac{1}{2}\theta})$ admits
a canonical action $\tilde\t\mapsto \tilde\beta_{\tilde\t}$ of $\tilde \Tb^N$ 
which permits the definition 
$$
C^\infty(M_{\theta},S):= \left( C^\infty(M,S)\widehat \otimes 
C^\infty(\Tb^N_{\frac{1}{2}\theta})
\right)^{\widehat\alpha\otimes\tilde\beta^{-1}},
$$
that is canonically a topological bimodule over $C^\infty(M_{\theta})$.
Next, since $[\Dslash , \widehat\alpha_{\tilde s}]=0$, the operator
   $$D_{\theta} := (\Dslash\otimes I) \downharpoonright C^\infty(M_{\theta},S)$$
is a  first-order differential operator on $C^\infty(M_{\theta},S)$.
There is also a hermitian structure 
$$
   (\psi\otimes t, \psi'\otimes t')= (\psi,\psi')\otimes t^\ast t' \ ,
$$
a $\mathbb Z_{2}$-grading $\gamma$ if $\mbox{dim}(M)$ is even and $M$ is oriented,
and an antilinear operator $\tilde J$, such that
$\tilde J(\psi\otimes t)=J\psi\otimes t^\ast$ for $\psi\in C^\infty(M,S)$,  
$t\in C^\infty(\tilde \Tb^N_{\theta})$. 

With these ingredients \cite{CD-V} define the spectral triple
$(C^\infty(M_{\theta}),\H_{\theta}, D_{\theta})$, where
$\H_\theta := \left(\H \widehat\otimes L^2(\tilde \Tb^N_{\theta})\right)
^{\widehat\alpha\otimes\tilde\beta^{-1}}$, and 
 $D_{\theta}$ is the closure of $\Dslash\otimes I$. Moreover,
$J_{\theta}$ and $\gamma$ are just (extensions of) previously-defined operators.\\
 ~\\
{\bf Proposition}~[CD-V].
The triple $(C^\infty(M_{\theta}),\H_{\theta}, D_{\theta})$ together
with the real structure $J_{\theta}$ and $\gamma$ satisfies all seven axioms 
required of a noncommutative manifold.\\
~\\
For instance, the summability is clear ($D_{\theta}$ is isospectral to $D$),
orientability comes from the invariance of $\vol_g$,
 which defines an invariant form on $M_{\theta}$
(an invariant Hochschild cycle in
$Z_{m}(A,A)$ for both $A=C^\infty(M)$ or $A=C^\infty(M_{\theta})$),
and Poincar\'e duality holds.\\

There is however a puzzling point about \eqref{spadesuit}.
Take for instance $M=\Tb^2$, $S=\Sigma_{0,0}$. 
According to \cite{CD-V} in view of \eqref{spadesuit},
 $$C^\infty(\Tb^2_{\theta},S)=
\left( C^\infty(\Tb^2,S)\widehat\otimes 
C^\infty(\Tb^2_{\theta/2})\right)^{\widehat\alpha\otimes\tilde\beta^{-1}}
$$
$$
=\left( C^\infty(\Tb^2)\otimes \bC^2 \widehat \otimes 
C^\infty(\Tb^2_{\theta/2})\right)^{\widehat\alpha\otimes{id}\otimes\tilde\beta^{-1}}
=C^\infty(\Tb^2_{\theta/2})\otimes \bC^2
$$
but on the other hand, according to `folklore wisdom'
$$
C^\infty(\Tb^2_{\theta},S)= C^\infty(\Tb^2_{\theta})\otimes \bC^2,
$$
(the trivial $C^\infty(\Tb^2_{\theta})$-module).
These two are evidently not isomorphic and
similar trouble happens for 
$\Tb^N$.
We resolve this puzzle by stating that 
\eqref{spadesuit} {\it does not} apply to the trivial spin structure $\sigma$ on tori.
This could be seen already from the immediate consequence of \eqref{spadesuit}
for $\te =0$ that 
\begin{equation}
\label{spadesuit0}
\tilde \Tb^N = \Tb^N .
\end{equation}
A glance at \eqref{cover} (and at the Examples) shows that \eqref{spadesuit0} selects 
a very special class of the covers which are 
twisted along one loop only.
This class corresponds to e.g. $\Sigma_{1,0}$ or $\Sigma_{0,1}$ on $\Tb^2$, 
more generally to $n$ spin structures (out of $2^N$) on $\Tb^N$, and to 
the two spin structures on $U(2)$ etc.
To know which cover it is concretely requires a supplementary information about the covering map. In any case \eqref{spadesuit0} excludes lot of spin structures, and even some manifolds like $M=S^N$. 
We shall see however that \eqref{spadesuit} for $\theta\neq 0$ put even stronger restrictions. Namely the matrix $\theta/2$ is possible as parameter matrix only for $N=2$ and otherwise it has to be more complicated.\\

Now we shall explain how and in which sense, when $N=2$, the `answer' \eqref{spadesuit} is indeed a deformation of such a special class of covers.
Moreover we will release the assumption
\eqref{spadesuit} and extend 
\cite{CD-V} to arbitrary spin structure 
by allowing any covering $\widetilde\Tb^2_{j_1,j_2}$ of 
$ \Tb^2$
~(and also generalize to $ \Tb^N$).
For this we really need to know what their $\theta$-deformations are.
Dually (using the concept of noncommutative principal $\bZ_2$-bundles),
$C^\infty((\widetilde\Tb^2_{j_1,j_2})_\te)$ and the embeddings of 
$C^\infty(\Tb^2_\te)$ as subalgebras of index 2 are given, respectively by:\\
$$
C^\infty(\Tb^2_\te) \otimes \bC^2 
\ni u_{1,\te}\otimes{_1\choose^1}, u_{2,\te}\otimes{_1\choose^1}
\leftarrow u_{1,\te},u_{2,\te}  ~{\rm if}~ j_1,j_2=0,0 \ ,  
$$
$$
C^\infty(\Tb^2_\frac{\te}{2}) 
\ni u^2_{1,\frac{\te}{2}},u_{2,\frac{\te}{2}}
\leftarrow u_{1,\te},u_{2,\te}  ~{\rm if}~ j_1,j_2=1,0 \ ,  
$$
$$
C^\infty(\Tb^2_\frac{\te}{2}) 
\ni u_{1,\frac{\te}{2}},u^2_{2,\frac{\te}{2}}
\leftarrow u_{1,\te},u_{2,\te}  ~{\rm if}~ j_1,j_2=0,1 \ ,  
$$
$$
C^\infty(\Tb^2_\frac{\te}{4})^{\bZ_2}
\ni u_{1,\frac{\te}{4}}^2, u_{2,\frac{\te}{4}}^2
\leftarrow u_{1,\te},u_{2,\te}  ~{\rm if}~ j_1,j_2=1,1 \ ,  
$$
where $\bZ_2$ acts by 
$
u_{1,\frac{\te}{4}}^m u_{2,\frac{\te}{4}}^n \mapsto
(-)^{m+n} u_{1,\frac{\te}{4}}^m u_{2,\frac{\te}{4}}^n 
$.
In fact, in view of \eqref{Z2prime},
$$ 
C^\infty(\Tb^2_\te) = C^\infty((\widetilde\Tb^2_{j_1,j_2})_\te)^{\bZ'_2}
$$
(the fixed elements), where the generator of $\bZ'_2$ acts by 
$$ 
\begin{array}{c}
a\otimes {_w\choose^z}\mapsto a\otimes {_z\choose^w}, ~{\rm if}~ j_1,j_2=0,0 \ , \\
\\
u^m_{1,\frac{\te}{2}} u^n_{2,\frac{\te}{2}}\mapsto 
(-)^m u^m_{1,\frac{\te}{2}} u^n_{2,\frac{\te}{2}},
~{\rm if}~ j_1,j_2=1,0 \ , \\
\\
u^m_{1,\frac{\te}{2}} u^n_{2,\frac{\te}{2}}\mapsto 
(-)^n u^m_{1,\frac{\te}{2}} u^n_{2,\frac{\te}{2}},
~{\rm if}~ j_1,j_2=0,1 \ , \\
\\
u_{1,\frac{\te}{4}}^m u_{2,\frac{\te}{4}}^n \mapsto
(-)^{mn} u_{1,\frac{\te}{4}}^m u_{2,\frac{\te}{4}}^n, 
~{\rm if}~ j_1,j_2=1,1\ 
\end{array}
$$
(in the last case $m+n$ is even).

Note that, although in the maximally twisted case the fourth root of 
$\lambda$ is involved in $ C^\infty(\Tb^2_\frac{\te}{4})$, 
only the square root matters in $C^\infty(\Tb^2_\frac{\te}{4})^{\bZ_2}$.
Anyway a kind of `transmutation' occurs: the more twisted $\sigma$ is, 
the more commutative parameter $\lambda$ is involved, 
namely $\lambda$, $\lambda^{1/2}$, $\lambda^{1/4}$. 
It is interesting to see if the sequence will continue for  $\Tb^N$, $N>2$ 
and in particular if it is the sequence 
$\lambda^{2^{-N}}\!\to\! 1$ as $N\!\to \!\infty$.

To answer this question consider $\Tb^N$ symmetries with $N\geq 3$.
Let $\te$ be a real skew $N\!\times\! N$ matrix of parameters. 
Using a straightforward generalization of \eqref{cover} and \eqref{Z2prime} 
to $N\geq 3$, the $\theta$-deformations of $\widetilde\Tb^N_{\underline j}$
have the following (dual) form, depending on the twist of $\sigma$
along the loops 
of $\Tb^N$.\\
If ${\underline j}={\underline 0}$ 
(the trivial cover), of course
$C^\infty(\widetilde\Tb^N_{\te}) 
=C^\infty(\Tb^N_{\te})\otimes\bC^2$.\\
If $j_m=1, j_i = 0$ for $i\neq m$, for some $1\leq m\leq N$
(twist only along the $m$-th loop), we have 
$C^\infty(\widetilde\Tb^N_{\te}) =  C^\infty(\Tb^N_{\widetilde\te})$, 
where
%
$$
\widetilde\te_{k,\ell}=
\left\{ \begin{array}{ll}
\te_{k,\ell}/2 ~&\mbox{if}~ k=m ~\mbox{or}~ \ell = m, \\
\te_{k,\ell} ~&\mbox{otherwise}
\end{array}
\right.
.
$$  
Note that $\widetilde\te \neq \te /2$ as a matrix 
(in contrast to the lowest case $N=2$).\\
If $j_i=1$ for $i\in X$ and $j_i = 0$ for $i\notin X$ 
(twist along $2\leq\! n\!\leq\! N$ loops with labels in $X$)
$C^\infty(\widetilde\Tb^N_{\te}) =  C^\infty(\Tb^N_{\widetilde\te})^{G_X}$, 
where 
$$
\widetilde\te_{k,\ell}=
\left\{ \begin{array}{ll}
\te_{k,\ell} ~&\mbox{if}~ k, \ell\notin X,\\
\te_{k,\ell}/2 ~&\mbox{if}~ k\in X, \ell\notin X ~\mbox{or vice versa},\\
\te_{k,\ell}/4 ~&\mbox{if}~ k,\ell\in X,
\end{array}
\right.
$$  
and 
\begin{equation}\label{GX}
G_X:= \{(\epsilon_1, \dots \epsilon_N)~|~ \epsilon_i = 1 ~\mbox{for}~ i\notin X,
\epsilon_i = \pm 1 ~\mbox{for}~ i\in X,~
\Pi_{i=1}^N\epsilon_i = 1\}
\end{equation}
($G_X $ is isomorphic to $(\bZ_2)^{n-1}\subset(\bZ_2)^n\subset(\bZ_2)^N $).
Concretely, $C^\infty(\Tb^N_{\tilde\te})^{G_X}$ consists of linear combinations
with rapidly decaying coefficients of the elements 
$u_{1,\tilde{\te}}^{k_1}u_{2,\tilde{\te}}^{k_2}\dots u_{N,\tilde{\te}}^{k_N}$
with $k_i+k_j$ even $\forall  i,j\in X$
(even in any pair of variables with labels in $X$). 
Notice that this is \underline{not} the even algebra in the variables with labels in $X$.
The embedding of $C^\infty(\Tb^N_{\te})$ is 
$u_{j,\te} \mapsto u^2_{j,\tilde{\te}}$ 
if $j\in X$ and 
$u_{j,\te} \mapsto u_{j,\tilde{\te}}$ 
if $j\notin X$.\\
~\\
Remark:
We see for arbitrary $\sigma$ (in particular for the maximally twisted one)
that at most the parameter $\te /4$ appears on the auxiliary level in 
$C^\infty(\Tb^N_{\widetilde\te})$.\\

Now everything works (with cumbersome but straightforward details) 
for arbitrary $\sigma$,
and we can state the following proposition.
%
%
\\
~\\
{\bf Proposition.} For arbitrary spin structure, there is a canonical action 
$\tilde\beta_{\tilde\t}$ of $\widetilde\Tb^N$ on $C^\infty(\widetilde\Tb^N_{\te}) $
(and on $L^2(\widetilde \Tb^N_{\theta}))$) and the spectral triple 
$$
\left(C^\infty(M_{\theta}),\left(\H \widehat\otimes L^2(\tilde \Tb^N_{\theta})\right)
^{\widehat\alpha\otimes\tilde\beta^{-1}}\!, D\otimes I\right)
$$ 
together with $\widehat J$ and $\gamma$ as before,
satisfies all axioms [C] of a noncommutative manifold.\\

We close with few remarks.\\
The spin structures $\Sigma_{j_1,j_2}$ on $\Tb^2$ and $\widetilde\Tb^2_{j_1,j_2}$ appear explicitly in \cite{D}.
\footnote{
There was no proofs correction and the long pre-publication errata list
was ignored by the publisher, who also introduced additional errors (e.g. page shifts). 
The merits include spinors on pseudo-Riemannian and on non-oriented manifolds, 
and the first (to the best of author's knowledge) definition of transformation 
of spinors under diffeomorphisms and equivariance 
$(\Dslash \psi)' = \Dslash' \psi'$ of the Dirac operator 
(from which invariance of the spectrum of $\Dslash$ is clear).}
\\
For spin structures on $C^\infty(\Tb^2_\te)$, dealt with in modern way, see \cite{PS}.\\
In this note, 
the interplay of $\te$-deformations and spin structures
was elucidated and,
filling the gaps, [CD-V, Sect.13] was extended to all spin structures $\sigma$ 
(e.g. to include $S^N$).\\
The extension of the method of [CL] via eigenmodes of $T^N$, $N>2$,
%
and of the (`splitting') isomorphism 
will be discused elswhere.\\
The study of spinor bundles for which the group $\bZ_2$ plays a key role,
can be generalized to bundles which admit action of $n$-fold coverings of tori.
This happens e.g. for $SU(n)$-isopinors, in which case the group $\bZ_n$ is relevant. 
As future work it will be interesting for the noncommutative geometry behind the standard model to study the case of $\bZ_3$,
and $\bZ_2 \times
\bZ_3
$~
(relevant e.g. for quarks which are color $SU(3)$-isospinors).

\end{document}